\documentstyle[12pt]{article}
\newcommand{\be}{\begin{equation}}
\newcommand{\ee}{\end{equation}}
\newcommand{\ben}{\begin{enumerate}}
\newcommand{\een}{\end{enumerate}}
\newcommand{\bit}{\begin{itemize}}
\newcommand{\eit}{\end{itemize}}
\newcommand{\ba}{\begin{array}}
\newcommand{\ea}{\end{array}}

\newtheorem{Thm}{Theorem}[section]
\newtheorem{Prop}[Thm]{Proposition}
\newtheorem{Lem}[Thm]{Lemma}
\newtheorem{Rem}[Thm]{Remark}
\newtheorem{Ex}[Thm]{Example}
\newtheorem{Def}[Thm]{Definition}
\newcommand{\dowod}{\noindent{\bf Proof:} }
\newcommand{\qed}{\hspace{5.5in} Q.E.D. }

\newcommand{\ov}{\overline}
\newcommand{\Ldwa}{\,\stackrel{2}{\bigwedge}}
\newcommand{\Ltrzy}{\,\stackrel{3}{\bigwedge}}
\newcommand{\w}{{\!}\wedge{\!}}
\newcommand{\Ph}{\mbox{\bf Ph}\, }
\newcommand{\End}{\mbox{\bf End}\, }
\font \msb=msbm10 scaled \magstep1
\font \mmsb=msbm8 scaled \magstep0
\font \dmsb=msbm10 scaled \magstep3

\newcommand{\rtimes}{\mbox{\msb o}\,}
\newcommand{\bR}{\mbox{\msb R} }
\newcommand{\bC}{\mbox{\msb C} }
\newcommand{\mbC}{\mbox{\mmsb C} }
\newcommand{\dbR}{\mbox{\dmsb R} }
\newcommand{\dbC}{\mbox{\dmsb C} }

\font \eul=eufm10 scaled \magstep2

\newcommand{\gotG}{\mbox{\eul g}}

\newcommand{\ar}{\alpha }

\newcommand{\br}{\beta }

\newcommand{\dr}{\delta }
\newcommand{\Dr}{\Delta }
\newcommand{\er}{\varepsilon }
\newcommand{\Lr}{\Lambda }
\newcommand{\lr}{\lambda }
\newcommand{\Om}{\Omega }

\newcommand{\sr}{\sigma }

\newcommand{\sgn}{{\rm sgn}\, }
\newcommand{\Hol}{{\em Hol}\, }

\newcommand{\cnabla}{\bar{\nabla}}
\newcommand{\cpartial}{\bar{\partial}}
\newcommand{\cz}{\bar{z}}

\begin{document}

\title{\bf Phase spaces related to standard classical $r$-matrices}
\author{{\bf S. Zakrzewski}  \\
\small{Department of Mathematical Methods in Physics,
University of Warsaw} \\ \small{Ho\.{z}a 74, 00-682 Warsaw, Poland} }

%\begin{table}[b]
%1991 Mathematical Subject Classification : $\;$ 22E46, 17B37
%\end{table}

\date{}
\maketitle

\begin{abstract}

Fundamental representations of real simple Poisson Lie groups
are Poisson actions with a suitable choice of the Poisson
structure on the underlying (real) vector space. We study these
(mostly quadratic) Poisson structures and corresponding phase
spaces (symplectic groupoids).

\end{abstract}

\section{Introduction}
The recent development of noncommutative geometry
and, in particular, the theory of quantum groups, raises the
question what happens with known models of physical systems when
we pass from usual configurations to non-commutative ones.
For the classical mechanical systems, this means that we allow
the configuration space to be a Poisson manifold (positions need
not commute). The phase space corresponding to usual
configuration manifold (Poisson structure equal zero) is its
cotangent bundle. For a general Poisson manifold, the
role of the phase space plays the corresponding symplectic
groupoid (if such exists, it is unique --- if one restricts to
connected and simply connected fibers).

It is natural to consider first mechanical systems with
symmetry. In the Poisson case a symmetry is described by a
Poisson action (of a Poisson group). This requirement imposes a
reasonable limitation on a choice of the Poisson structure and
actually leads to a construction of it.

In this paper we construct Poisson structures on  real
finite-dimensional vector spaces (the configuration spaces),
such that the action of a chosen linear simple Poisson group
becomes a Poisson action (the Poisson structure on the group is
typically given by a standard classical $r$-matrix). We
construct also the corresponding phase spaces.

\section{Preliminaries and notation}
For the theory of Poisson Lie groups we refer to
\cite{D:ham,D,S-T-S,Lu-We,Lu}. Let us recall some basic notions
and facts. We follow the notation used in our previous papers
\cite{repoi,poihom,PPgr}.

A {\em Poisson Lie group} is a Lie group $G$ equipped with a
Poisson structure $\pi $ such that the multiplication map is
Poisson. The latter property is equivalent to the following
property (called {\em multiplicativity} of $\pi$)
$$
\pi (gh)=\pi (g)h+g\pi (h) \qquad \mbox{for}\;\; g,h\in G.
$$
Here $\pi (g)h$ denotes the right translation of $\pi (g)$ by
$h$ etc. This notation will be used throughout the paper.

A Poisson Lie group is said to be {\em coboundary} if
\be\label{cob}
\pi (g) = rg -gr
\ee
for a certain element $r\in \Ldwa \gotG$. Here $\gotG$ denotes
the Lie algebra of $G$. Any bivector field of
the form (\ref{cob}) is multiplicative. It is Poisson if and
only if
$$
[r,r]\in (\Ltrzy \gotG )_{\rm inv}
$$
(the Schouten bracket $[r,r]$ is $\gotG$-invariant).
In this case the element $r$ is said to be a {\em classical
$r$-matrix} ({\em on} $\gotG$).

If $G$ is semisimple, any Poisson Lie group structure on $G$ is
coboundary. Standard classical $r$-matrix for a simple group ---
such that corresponds to `the standard (quantum)
$q$-deformation' --- is given by (cf. \cite{BD}, Prop. 2.1 in
\cite{J})
\be\label{stand}
r = c\sum_{\ar >0} \frac{X_{\ar}\w X_{-\ar}}
{\left\langle X_{\ar },X_{-\ar} \right\rangle },
\ee
where $X_{\pm \ar}$ are (positive and negative) root vectors
relative to a Cartan subalgebra in $\gotG$,
$\left\langle \cdot ,\cdot \right\rangle $ is the Killing form
and $c$ is a constant (if $G$ is compact, $X_{-\ar } = \ov{X}
_\ar$ and $c$ is imaginary).

Let $(G,\pi )$ be a Poisson Lie group. An action of $G$ on a Poisson
manifold $(M, \pi _M)$ is said to be a {\em Poisson action} if
the action map $G\times M\to M$ is Poisson. It holds if and only
if the following $(G,\pi)$-{\em multiplicativity} of $\pi _M$ is
satisfied
$$
\pi_M (gx) = \pi (g)x+g\pi _M (x)\qquad    \mbox{for}\;\; g\in
G,x\in M.
$$
For any fixed action $G\times M\ni (g,x)\mapsto gx\in M$ and any
$k$-vector $w\in {\,\stackrel{k}{\bigwedge}}\, \gotG $ we denote
by $w_M$ the associated $k$-vector field on $M$:
$$
w_M (x):= wx.
$$

\section{Problem}\label{problem}
The classical $r$-matrices for simple Lie groups like $SL(n,\bR
)$, $SO(n,\bR )$, $SU(n)$ are relatively well investigated (in
the sequel we shall consider mainly the standard $r$-matrices
(\ref{stand}), which are indeed representing the non-trivial part
of all classical $r$-matrices). In order to consider mechanical
systems based on Poisson symmetry (typically being a `deformation'
of some ordinary symmetry), we have first to deal with
the following problems:
\ben
\item Given an action $G\times M\to M$ (the ordinary symmetry) and
a Poisson structure $\pi $ on $G$ making it a Poisson Lie group
$(G,\pi )$ (a `deformation' of the group), find all Poisson
structures $\pi _M$ on $M$ such that the  action becomes Poisson (the
`deformed' symmetry).
\item In cases when $M$ plays the role of the configurational
manifold, construct the phase space $\Ph (M,\pi _M)$ i.e. the
symplectic groupoid of $(M,\pi _M)$.
\een
For symplectic groupoids, phase spaces of Poisson manifolds etc.
we refer to \cite{CDW,Ka,qcp,We,srni,X}.

For simplicity, in this paper we consider only the essential
part of the structure of the symplectic groupoid (which is
in most cases sufficient to formulate the classical model). Namely,
for a given Poisson manifold $(M,\pi _M)$ of dimension $k$ we shall
construct a symplectic manifold $S$ of dimension $2k$,
a surjective Poisson map from $S$ to $M$ and a lagrangian
section of it. In this case, we shall simply call $S$ the {\em
phase space} of $(M,\pi _M)$.

\section{Fundamental bi-vector field}

Let $G\times M\to M$ be an action. Let $r\in \Ldwa \gotG$ be
a classical $r$-matrix and $\pi$ --- the corresponding Poisson
structure (\ref{cob}) (this notation is fixed throughout the
section).
\begin{Lem}\label{L1}
\ben
\item $r_M$ \ \  is $(G,\pi )$-multiplicative
\item any $(G,\pi )$-multiplicative $\pi _M$ is given by
$\pi _M = r_M + \pi _{\rm inv}$, where $\pi _{\rm inv}$  is a
$G$-invariant bi-vector field
\item\label{L13} $[r_M,\pi _{\rm inv}]=0$
\item\label{rM} $[r_M,r_M]=[r,r]_M$.
\een
\end{Lem}
Point 1 follows from $r(gx) = (rg -gr)x + g(rx)$. Point 3
follows from the fact that $r_M$ is built out of the fundamental
vector fields of the action (and these vector fields preserve
$\pi _{\rm inv}$). From 3 follows that if both $r_M$ and $\pi
_{\rm inv}$ are Poisson then $\pi _M$ is also Poisson.
Point 4 follows from the known property of fundamental fields of
the action:
$$
[X_M,Y_M] = [X,Y]_M\qquad\qquad\mbox{for} \;\; X,Y\in \gotG
$$
(the Lie bracket on $\gotG$ being defined by identifying
elements of $\gotG$ with the corresponding {\em right-invariant}
vector fields on $G$).

In analogy with fundamental vector fields $X_M$, we call $r_M$
the {\em fundamental bi-vector field}. It is essential to know
whether it is Poisson.

\begin{Ex}{\rm {\sl Poisson Minkowski spaces (Poincar\'{e} group
action)}. \
Any invariant element of $\Ltrzy \gotG$, where $\gotG =\bR
^4\rtimes o(1,3)$ is the Poincar\'{e} Lie algebra, is
proportional to
\be\label{Om}
\Om = g ^{jk} g ^{lm} e_j\wedge e_l \wedge \Om _{km} \qquad
\qquad \Om _{km}:= e_k\otimes g(e_m) - e_m\otimes g(e_k)\in o(1,3),
\ee
(summation convention), where $(e_j)_{j=0,\ldots,3}$ is a basis
in $M=\bR ^4$, $g$ is the Lorentz metric and $g^{jk}$ are the
components of the contravariant metric (cf. \cite{PPgr,PsPg}). Since
$$
\Om _M (x)= g ^{jk} g ^{lm} e_j\wedge e_l \wedge
(e_kg(e_m,x)-e_mg(e_k,x))=0,
$$
for each classical $r$-matrix on $\gotG$ the fundamental bi-vector field
$r_M$ on $M$ is Poisson (because $[r_M,r_M]=[r,r]_M\sim \Om
_M=0$). By point 2 of Lemma~\ref{L1} this is the only $(G,\pi
)$-multiplicative bivector field on $M$,
since zero is the only $G$-invariant bivector field on $M$.
(Recall also that any Poisson structure on $G$ comes from an
$r$-matrix \cite{PPgr}.) Concluding: for each Poisson
Poincar\'{e} group there is exactly one Poisson Minkowski space
(see also \cite{poihom}). All this is true also for the case of
arbitrary signature, $\gotG =\bR ^{p+q}\rtimes o(p,q)$, in
dimension $n=p+q>3$. (Cf.~\cite{PW} for the quantum case.)  }
\end{Ex}

\begin{Ex}{\rm {\sl Poisson Minkowski spaces (Lorentz group
action)}. \
Classical $r$-matrices for the Lorentz Lie algebra $o(1,3)$ are
classified in \cite{repoi}. We know that $[r,r]=[r_-,r_-]$ and it
is non-zero only in the case $r_-=i\lr X_+\wedge X_-$ (in the
classification of \cite{repoi}) with $\lr \neq 0$:
$$
[r_-,r_- ]= -\lr ^2 [X_+\wedge X_- ,X_+\wedge X_-]= 2\lr ^2
X_+\wedge [X_+,X_-]\wedge X_- =
4\lr ^2
X_+\wedge H\wedge X_- ,
$$
where  $X_+$, $X_-$, $H$ is the standard basis:
$$H=\frac12 \left[ \ba {cc} 1 & 0 \\ 0 & -1 \ea \right] ,\qquad
X_+=\left[ \ba {cc} 0 & 1 \\ 0 & 0 \ea \right] ,\qquad
X_-=\left[ \ba {cc} 0 & 0 \\ 1 & 0 \ea \right] .$$
Considering the usual action of the Lorentz Lie algebra on the
Minkowski space $M=\bR ^{1+3}$, we obtain
$$ (X_+\wedge X_- )_M(x) = 2\Om _{01} (x)\wedge\Om _{13}(x) ,$$
where (see (\ref{Om}))
$$ \Om _{jk} (x) = e_j x_k - e_k x_j.$$
Since  $\Om _{jk}(x)$, $\Om _{kl}(x)$ and $\Om _{lj}(x)$
are linearly dependent for each fixed $j,k,l$,
$$
(X_+\wedge H\wedge X_- )_M(x) =-2\Om _{30}(x)\wedge \Om _{01}(x)
\wedge\Om _{13}(x)=0
$$
($H_M(x)= \Om _{30}(x)$), but
$$
(X_+\wedge JH\wedge X_- )_M(x) =-2\Om _{21}(x)\wedge \Om _{01}(x)
\wedge\Om _{13}(x)
$$
($J$ is the complex structure in $\gotG$) is not zero. It
follows that $r_M$ is Poisson if and only if
$\lr ^2$ is real, i.e. either $\ar$ or $\br$ in \cite{repoi} has
to be zero. Moreover, since the only Lorentz invariant bivector
field on $M$ is zero, $r_M$ is the only $(G,\pi
)$-multiplicative field on $M$. It follows that for
$\ar\cdot\br\neq 0$ there is no Poisson structure on $M$ such
that the action is Poisson. (Similar fact should hold for quantum
Lorentz groups \cite{cl}: $q$ should be real or of modulus one.)
}
\end{Ex}

Returning to a general technique,
consider now two special cases of $r$-matrices.

\subsection{Triangular case: $[r,r]=0$}           \label{triang}

 Let $\xi \colon T^*M\to
M$ be the cotangent bundle projection and let $\pi _0$ denote
the canonical Poisson structure of $T^*M$. In the triangular case
  \ben
   \item $r_M$ is Poisson (by Lemma~\ref{L1}.\ref{rM})
   \item $r_{T^*M}$ is Poisson (also Lemma~\ref{L1}.\ref{rM}); \ $\xi  _*
r_{T^*M}=r_M$
   \item $\pi _{T^*M}:=r_{T^*M}+\pi _0$ is Poisson (by
Lemma~\ref{L1}.\ref{L13}); \ $\xi  _*\pi _{T^*M}=r_M$.
  \een
This means that problems formulated in Sect.~\ref{problem} are
relatively easily solved. As the phase space one can take the
open subset of points in $T^*M$, in which the Poisson structure
$\pi _{T*M}$ is non-degenerate (it is certainly non-degenerate
in a neighbourhood of the zero section --- that's why we have
added $\pi _0$ in 3). (To construct the
symplectic groupoid one should still find the foliation symplectically
orthogonal to the fibers of the projection and choose points
which have also the projection on $M$ along this foliation.)

For another approach to this case, see \cite{X}.

 \subsection{Case of a simple $\gotG$}

 In this case one can use
the method of \cite{GP} to rewrite the condition $[r_M,r_M]=0$.
Denote by $\Omega $ the canonical invariant element of $\Ltrzy
\gotG$. Its Killing transported version to $\Ltrzy \gotG^*$ is
defined by
$$
\Omega ^\dagger (X,Y,Z)= <[X,Y],Z>.
$$
It is known \cite{ccc} that all invariant elements of $\Ltrzy\gotG
$ are proportional to $\Omega$, hence $[r,r]\sim \Omega$.
Suppose $[r,r]$ is not zero. Then $r_M$ is Poisson
$\Leftrightarrow$ $\Omega _M =0$
(in general, $\Omega _M$ is just $G$-invariant).
Now,

$\Omega x=0 \Leftrightarrow $ the composition of linear maps
$$ \bR \stackrel{\Omega}{\to} \Ltrzy \gotG \to \Ltrzy (\gotG /
\gotG _x)$$
is zero $\Leftrightarrow$
the composition of linear maps
$$ \bR \stackrel{\Omega ^*}{\leftarrow} \Ltrzy \gotG \leftarrow\Ltrzy
(\gotG _x)^0$$
is zero $\Leftrightarrow$
the composition of linear maps
$$ \bR \stackrel{\Omega ^\dagger }{\leftarrow} \Ltrzy \gotG
\leftarrow \Ltrzy  (\gotG _x )^\perp $$
is zero $\Leftrightarrow \Omega ^\dagger \vert _{(\gotG
_x)^\perp } =0\Leftrightarrow  <[X,Y],Z>=0$ for $X,Y,Z\in \gotG
_x^\perp \Leftrightarrow [\gotG _x^\perp ,\gotG _x^\perp
]\subset \gotG _x $.
Concluding,
\be\label{crit}
[r_M,r_M]x=0\Longleftrightarrow
[\gotG _x^\perp ,\gotG _x^\perp ]\subset \gotG _x .
\ee
The advantage of this method is that we do not have to use the
explicit form of $\Om $.
\begin{Prop}\label{class}
In the following three cases, for any classical $r$-matrix on
$\gotG$ the fundamental bi-vector field $r_M$ on $M$ is Poisson:

1. $\gotG = so(n,\bR )$, $M=\bR ^n$ \hspace{5mm}
2. $\gotG = sl(n,\bR )$, $M=\bR ^n$ \hspace{5mm}
3. $\gotG = sp(n,\bR )$, $M=\bR ^{2n}$ \hspace{5mm}.

\noindent
For  $\gotG = su(n)$, $M=\bC ^n=\bR ^{2n}$, the fundamental
bi-vector field is Poisson if and only if $r$ is triangular.
\end{Prop}

\dowod Let $e_1,\ldots ,e_n$ be the standard basis in $\bR^n$.
The dual basis in $(\bR ^n )^*$ is denoted by $e^1,\ldots ,e^n$.
We consider subsequent cases separately.

1. $\gotG = so(n,\bR )$. Set $x=e_n$, then $\gotG _x$ is spanned by
$\Om _{jk}$ for
$j,k<n$ and $\gotG _x^{\perp}$ is spanned by $\Om _{jn}$ for
$j<n$. It is easy to see that $[\Om _{jn},\Om _{kn}]$ is
proportional to $\Om _{jk}$, hence it belongs to $\gotG _x$.
This shows (by (\ref{crit})) that $[r_M,r_M]=0$.

2. $\gotG = sl(n,\bR )$. For the same $x=e_n$, $\gotG _x$
 is spanned by
$e_n{^j}$ (the matrix units) for $j<n$, but
$[e_n{^j},e_n{^k}]=0$.

3. $\gotG = sp(n,\bR )$. We use the basis $e_1,\ldots
,e_n,e^1,\ldots ,e^n$ in $M=\bR
^n\oplus (\bR ^n)^* \cong \bR ^{2n}$. We have four types of matrix units
defined by
$$
e_j{^k}:= e_j\otimes e^k,\qquad e^j{_k}:=e^j\otimes e_k,\qquad
e_{jk} := e_j\otimes e_k,\qquad e^{jk}:= e^j\otimes e^k,
$$
with the action on $(x,p)\in M=\bR ^n\oplus (\bR ^n)^*$ given
explicitly by
$$
e_j{^k} (x,p)=e_jx^k,\qquad e^j{_k}(x,p) = e^jp_k,\qquad
e_{jk} (x,p) = e_jp_k,\qquad e^{jk} (x,p)=e^jx^k.
$$
We use the following basis in $\gotG =sp(n,\bR )$:
\be\label{basis}
a_{jk}:= e_{jk}+e_{kj} \;\;(j\leq k),\qquad b^{jk}:=
e^{jk}+e^{kj} \;\; (j\leq k),\qquad d_j{^k}:= e_j{^k} - e^k{_j}.
\ee
For $x:= e_n$, $\gotG _x$ is spanned by $a_{jk}$, $b^{jk}$ with
$j,k<n$ and $d_j{^k}$ with $k<n$ and $\gotG _x ^{\perp }$ is
spanned by $a_{jn}$, $d_n{^j}$ with $j=1,\ldots n$.
Now, $[a_{jn}, a_{kn}]=0$ and $[d_n{^j},a_{kn}]=\dr^j_ka_{nn} +
\dr ^j_n a_{kn}\in \gotG _x$.

We now pass to the case of $\gotG = su(n)$, with the basis
\be\label{sunbas}
F_j{^k} := e_j{^k} - e_k {^j},\qquad G_j{^k} := i(e_j{^k} + e_k
{^j}),\qquad H_j := e_{j+1}{^{j+1}} - e_j{^j}.
\ee
For $x=e_n$, $\gotG _x$ is spanned by $F_j{^k}$, $G_j{^k}$ for
$j,k<n$ and $H_j$ for $j<n-1$ and $\gotG _x ^{\perp
}$ contains $F_j{^n}$, $G_j{^n}$ for $j<n$ (and the component of
$H_{n-1}$, orthogonal to $H_j$ for $j<n-1$).
We have
$$
[F_1{^n} ,G_1{^n}]=2i (e_1{^1}-e_n{^n})\not\in \gotG _x .
$$
Actually one can see that
$$
[\gotG _x^{\perp} ,\gotG _x ^{\perp}]= \gotG _x + \left\langle
i(e_1{^1}-e_n{^n})\right\rangle
$$
(inclusion (\ref{crit}) is violated only by one dimension).

\qed

\section{Standard $r$-matrices for $sl(n,\dbR )$ and $sp(n,\dbR )$}

According to (\ref{stand}),
the standard $r$-matrix for $\gotG = sl(n,\bR )$ is given by
\be\label{stasln}
r = \er\sum_{j<k} e_j{^k}\wedge e_k{^j} \qquad\qquad (\er\in\bR )
\ee
(the Cartan subalgebra consists of diagonal matrices and the
`positive' roots are contained in the upper-triangular
matrices). Considering the natural action of $\gotG$ on $M:=\bR
^n$, we obtain
\be\label{slrMx}
r_M(x) =\er \sum_{j<k}x^jx^k e_j\wedge e_k = \sum_{j<k}(x^j
e_j)\wedge (x^ke_k),
\ee
which defines the following Poisson brackets of coordinates:
\be\label{slxx}
\{ x^j,x^k\} =\er x^jx^k\qquad\qquad \mbox{for}\;\; j<k .
\ee

For $\widetilde{\gotG}:=sp(n,\bR)$ we choose $d_j{^k}$ $(j<k)$
and $a_{jk}$ $(j\leq k)$ as positive roots, which gives the
following expression for the standard $r$-matrix (we denote it
by $\widetilde{r}$):
$$
\widetilde{r} = \er\left(\sum_{j<k} d_j{^k}\wedge d_k{^j}
+ \frac{1}{2}\sum_{j,k} a_{jk}\wedge b^{jk} \right)
$$
(notation as in (\ref{basis})). Considering the natural action
of $\widetilde{\gotG}$ on $\widetilde{M}=T^*M=\bR ^n\oplus (\bR
^n)^*$, we obtain
\be\label{wr}
\widetilde{r}_{\widetilde{M}}(x,p) =\er
\left[\sum_{j<k} (e_jx^k-e^kp_j)\wedge (e_kx^j - e^jp_k)
+ \frac{1}{2}\sum_{j,k} (e_jp_k + e_kp_j)\wedge (e^jx^k+e^kx^j)
\right]
\ee
$$
= \er\left[x\wedge p + \sum_{j<k} (x^jx^ke_j\wedge e_k -
p_jp_ke^j\wedge e^k) +\sum_{j} \left(\sum_{k} (1 -\sgn (k-j) )
x^kp_k\right) e_j\wedge e^j \right],
$$
which gives the following quadratic Poisson brackets of
coordinates and momenta
\be\label{xxpp}
\{ x^j,x^k\} =\er x^jx^k,\qquad \{p_j,p_k\} = -\er
p_jp_k\qquad\mbox{for}\;\; j<k
\ee
\be\label{mixed}
\{ x^j,p_k\} = \er\left[x^jp_k + \dr ^j{_k} \left(\sum_{i} (1 -\sgn (i-j) )
x^ip_i\right)\right] .
\ee
Now observe that we can add to (\ref{wr}) the canonical
bi-vector $\pi _0$ on $T^*M$ which modifies (\ref{mixed}) in the
following way
\be\label{mixed1}
\{ x^j,p_k\} = -\dr ^j{_k}+\er\left[x^jp_k + \dr ^j{_k}
\left(\sum_{i} (1 -\sgn (i-j) ) x^ip_i\right)\right] .
\ee
Poisson structure (\ref{xxpp}), (\ref{mixed1}) projects on
(\ref{slxx}) and is non-degenerate
in a neighbourhood of $(\bR ^n \oplus \{0\})\cup (\{0\}\oplus
(\bR ^n)^*)$. We have thus constructed
the phase space for $(M,r_M )$ ($M$ is indeed embedded into this
phase space as a lagrangian submanifold: $p_j=0$
are {\em first class} constraints).

The quantum version of the above construction has been
described in \cite{Zu}.

\begin{Rem}
{\rm The natural embedding of $\gotG$ into $\widetilde{\gotG}$
(the lift to $T^*M$),
given by $e_j{^k}\mapsto d_j{^k}$, is a homomorphism of Lie
bialgebras. It follows that the action of $SL(n,\bR )$ on $T^*M$
is Poisson.}
\end{Rem}

\begin{Rem}  {\rm
Set $\nabla _j:= x^je_j $ (no summation). According to
(\ref{slrMx}), $r_M = \sum_{j<k} \nabla _j\wedge \nabla _k$. We
see that $r_M$ is built of commuting vector fields and, actually,
$r_M=\rho _M$, where $\rho $ is a $r$-matrix on the abelian Lie
algebra spanned by those fields. Since $\rho$ is triangular and
even abelian \cite{abel}, one can easily construct the phase
space of $(M,\rho _M)$ using the method described in
Sect.~\ref{triang}, or, another method, described in
\cite{abel}. In the present paper we shall exploit only the
original $r$-matrix, because such an approach can be generalized
and applied to more difficult situations (see next sections).}
\end{Rem}

\section{Crossed product phase spaces and quasitriangularity}

Let $\gotG$ be a Lie subalgebra of $\End V$, where $V=\bR ^n$.
Any classical $r$-matrix on $\gotG$ may be identified as an
element of $\Ldwa \End V $, which can be expressed in terms of
matrix units:
$$
 r = \sum_{jklm} r^{jk}_{lm} e_j{^l}\otimes e_k{^m}\qquad\qquad
r^{jk}_{lm} = -r^{kj}_{ml}.
$$
If $r$ is triangular, the phase space of $(V,r_V)$ can be
realized on $T^*V=V\times V^*$ with the Poisson structure $\pi
_{T^*V}$ being the sum of the canonical Poisson structure $\pi
_0$ and $r_{T^*V}$ (cf.~Sect.~\ref{triang}). We have then
$$
\pi _{T^*V} (x,p)=  \pi _0 + \sum_{jklm} r^{jk}_{lm}
(e_jx^l-p_je^l)\otimes (e_kx^m-p_ke^m),
$$
which leads to the following Poisson brackets
\be\label{PB1}
\{ x^j,x^k\} = \sum_{lm} r^{jk}_{lm}x^lx^m,\qquad \{ p_l,p_m\} =
\sum_{jk} p_jp_kr^{jk}_{lm}
\ee
\be\label{PB2}
\{ x^k,p_l\} = -\dr ^k{_l}+\sum_{jm}p_j r^{jk}_{lm}x^m .
\ee
We shall use also the following abbreviated notation:
\be\label{PBt}
\{ x_1,x_2\} = rx_1x_2,\qquad \{ p_1,p_2\} = p_1p_2r,\qquad
\{x_1,p_2\} = -I + p_1rx_2.
\ee
\begin{Def} {\rm  Let $(M,\pi _M )$ and
$(N,\pi _N )$ be two Poisson manifolds and $P:= M\times N$. If
$\pi _P$ is a Poisson structure on $P$ such that the cartesian
projections $P\to M,N$ are Poisson, then $(P,\pi _P )$ is said
to be a} crossed product {\rm of $(M,\pi _M )$ and $(N,\pi _N )$.}
\end{Def}
We see that in the triangular case, the phase space is a crossed
product of $(V,r_V )$ and $(V^*,r_{V^*})$. Note that the Poisson
brackets between $x^k$ and $p_l$ (the cross-relations) are expressed in
terms of the same $r$-matrix as $r_V$ and $r_{V^*}$.

In the case of the non-triangular $r$-matrix (\ref{stasln}) for
$sl(n,\bR )$, the phase space turns out to be also a crossed
product of $(V,r_V )$ and $(V^*,r_{V^*})$. Poisson structure
(\ref{xxpp}), (\ref{mixed1}) is realized on $T^*V$ and has the
following form
\be\label{+Delta}
\pi _{T^*V}=\pi_0 + r_{T^*V} + \Dr ,
\ee
where $\Dr $ is some additional quadratic term (in the
cross-relations), which was not present in the triangular case.

We shall now explain the nature of the additional term $\Dr$
in a general situation. We assume that our $r$-matrix is
such that $r_V$ and $r_{V^*}$ are Poisson (with explicit form of
Poisson brackets given by (\ref{PB1})) and we ask what
conditions should satisfy a bi-vector field $\Dr $ of the form
\be\label{De}
\Dr (x,p) = \sum_{jklm} p_j\Dr ^{jk}_{lm}x^m e_k\wedge e^l
\ee
on $T^*V$, in order to make (\ref{+Delta}) a Poisson bi-vector field.
Since $r_{T^*V}(x,p)= r_V(x)+r_{V^*}(p) +\sum_{jklm} p_jr
^{jk}_{lm}x^m e_k\wedge e^l$, we have
$$
\pi _{T^*V}(x,p)=\pi_0 + r_{V} +r_{V^*} + \sum_{jklm} p_j w
^{jk}_{lm}x^m e_k\wedge e^l,
$$
where $w^{jk}_{lm} = r^{jk}_{lm} + \Dr ^{jk}_{lm}$.
We look therefore for conditions on $w^{jk}_{lm}$ under which
the brackets
\be\label{general}
\{ x_1,x_2\} = rx_1x_2,\qquad \{ p_1,p_2\} = p_1p_2 r ,\qquad
\{p_1,x_2\} = I - p_1wx_2
\ee
(abbreviated notation) satisfy the Jacobi identity. We first
consider only the quadratic part
\be\label{quadr}
\{ x_1,x_2\} = rx_1x_2,\qquad \{ p_1,p_2\} = p_1p_2 r ,\qquad
\{p_1,x_2\} =- p_1wx_2
\ee
(it is easy to show that the quadratic part itself must
also  define a Poisson bracket).
Since
$$
\{p_1,\{ x_2,x_3\}\} = \{ p_1,r_{23}x_2x_3\} = - p_1
r_{23}(w_{12}+w_{13})x_1x_3
$$
and
$$
\{ \{p_1,x_2\} ,x_3\} - \{\{p_1,x_3\},x_2\} =
\{ -p_1w_{12}x_2 ,x_3\} - \{ p_1w_{13}x_3 ,x_2\}
$$
$$
= p_1 (w_{13}w_{12} - w_{12}r_{23})x_2x_3 -
(p_1 (w_{12}w_{13} - w_{13}r_{32})x_2x_3),
$$
the part of the Jacobi identity corresponding to the equality $\{p_1,\{
x_2,x_3\}\} = \{ \{p_1,x_2\} ,x_3\} - \{\{p_1,x_3\},x_2\} $ is
equivalent to
\be\label{Jac1}
p_1 ([w_{12},w_{13}] + [w_{12}, r_{23}] + [w_{13}, r_{23}])
x_2x_3 =0.
\ee
Similarly, the part of the Jacobi identity corresponding to $\{ \{
p_1,p_2\},x_3\} = \{ p_1, \{ p_2,x_3\}\}-\{ p_2,\{p_1,x_3\}\}$
is equivalent to
\be\label{Jac2}
p_1p_2 ([r _{12},w_{13}] + [r _{12}, w_{23}] + [w_{13},
w_{23}])x_3=0.
\ee
\begin{Thm}\label{tw}
Let $(\gotG , r)$ be quasitriangular, i.e. there exists an
invariant symmetric element $s$ of $\gotG \otimes\gotG$ such
that $w:= r + s$ satisfies the classical Yang-Baxter equation:
$$
[[ w, w]] := [w_{12},w_{13}] + [w_{12}, w_{23}] + [w_{13}, w_{23}]=0.
$$
Then the brackets (\ref{quadr}) satisfy the
Jacobi identity (we assume that (\ref{PB1}) satisfy already the
Jacobi identity).
\end{Thm}
\dowod
 Since $s$ is invariant,
$$
[w_{12}, s_{23}] + [w_{13},s_{23}]=0
$$
 (for any $w$), hence  $[[w,w]]=0$ if and only if
\be\label{subtr}
[w_{12},w_{13}] + [w_{12}, r_{23}] + [w_{13}, r_{23}]=0.
\ee
This obviously implies (\ref{Jac1}). Similarly, since $
[s_{12},w_{13}] + [s_{12}, w_{23}]  =0$, $[[w,w]]$ is zero if and
only if
$$
[r_{12},w_{13}] + [r_{12}, w_{23}] + [w_{13},w_{23}] =0,
$$
which implies (\ref{Jac2}).

\qed
\begin{Rem} {\rm  Let $r$, $s$ and $w$ satisfy the assumptions
of Theorem~\ref{tw}.  For any $\lr \in \bR$, the element
\be
w_{\lr } := r+s+\lr I\otimes I\in \End V\otimes \End V,
\ee
satisfies the Yang-Baxter equation and the proof of
Theorem~\ref{tw} works also for $w_\lr$ (once $w$
satisfies (\ref{subtr}),  $w_{\lr}$
also satisfies (\ref{subtr})).
It means that one can replace
$w$ by $w_\lr $ in (\ref{quadr}).  }
\end{Rem}
\begin{Thm}
Under assumptions of Theorem~\ref{tw}, brackets (\ref{general})
satisfy the Jacobi identity if and only if
\be\label{syms}
s^{jk}_{lm} = s^{kj}_{lm}.
\ee
Brackets (\ref{general}) with $w$ replaced by $w_\lr $ satisfy
the Jacobi identity if and only if
\be\label{symsl}
(s_\lr )^{jk}_{lm}  = (s_\lr )^{kj}_{lm},
\ee
where $s_\lr = s +\lr I\otimes I$.
\end{Thm}
\dowod In
the part of the Jacobi identity corresponding to $\{p_1,\{
x_2,x_3\}\} = \{ \{p_1,x_2\} ,x_3\} - \{\{p_1,x_3\},x_2\} $, we
must take care of the linear terms (the cubic terms were taken
into account in the previous theorem). This gives
\be\label{r-w}
rx - (rx)^{\rm t} = wx - (wx)^{\rm t},
\ee
where $(rx)^{jk}_l := \sum_m r^{jk}_{lm}x^m$, $((rx)^{\rm
t})^{jk}_l:= (rx)^{kj}_l$, etc. Of course, (\ref{r-w}) means that
$sx = (sx)^{\rm t}$, i.e. (\ref{syms}). The
modification $w\mapsto w+\lr I\otimes I$, $s\mapsto s+\lr
I\otimes I$ leading to condition (\ref{symsl}) is straightforward.
The remaining part of the Jacobi identity leads to the same condition.

\qed
\begin{Ex}
{\rm It is convenient to consider (\ref{stasln})
as a $r$-matrix on $\gotG =gl (n,\bR )$, because it is easy to write
down a natural invariant symmetric element (trace form) of
$\gotG\otimes \gotG$. Taking this element with the appropriate coefficient:
$$
 s= \er \sum_{j,k} e_j{^k}\otimes e_k{^j},
$$
we obtain
$$
  w=r+s=\er \left( \sum_{j} e_j{^j}\otimes e_j{^j} + 2\sum_{j<k}
e_j{^k}\otimes e_k{^j} \right) ,
$$
which satisfies the classical Yang-Baxter equation.
There is a unique modification $s_\lr =s+\lr I\otimes I$ of $s$
satisfying the symmetry (\ref{symsl}), namely for $\lr =\er$:
$$
 s_\lr = s_\er = s +\er I\otimes I,\qquad
(s_\er)^{jk}_{lm}=\er(\dr ^j_m\dr ^k_l + \dr ^j_l\dr ^k_m ).
$$
Poisson brackets (\ref{general}) with $w_\lr = w_\er $ coincide with
(\ref{xxpp}), (\ref{mixed1}). }
\end{Ex}

Concluding, if $(\gotG ,r)$ is quasitriangular and if there
exists a modification  $s_\lr =s+\lr I\otimes I$ of $s$
satisfying (\ref{symsl}), then one can realize the phase space
of $(V,r_V)$ on $T^*V$ with the Poisson structure
(\ref{+Delta}), where $\Dr $ is given by (\ref{De}) with
\be
\Dr ^{jk}_{lm} = (s_\lr ) ^{jk}_{lm}.
\ee
It is convenient to introduce the following notation. For each
\be
\rho = \sum_{jklm} \rho ^{jk}_{lm} e_j{^l}\otimes e_k{^m}\in
\End V\otimes \End V,
\ee
we denote by $\rho _{VV^*}$ the bi-vector field on $T^*V=V\oplus
V^*$ defined by
\be
\rho _{VV^*} (x,p) = \sum_{jklm} p_j\rho ^{jk}_{lm}x^m e_k\wedge
e^l.
\ee
Using this notation we can write (\ref{+Delta}) as follows
\be
\pi _{T^*V} = \pi _0 + r_{T^*V} + (s_\lr )_{VV^*}=\pi _0 +
r_V+r_{V^*} + (w_\lr ) _{VV^*}.
\ee

\section{$SO(n,\dbR )$, imaginary quasitriangularity and reality
condition}

In this section we construct the phase space of $(V, r_V)$,
where $V:=\bR ^n$, for a standard $r$-matrix on $so(n,\bR )$,
using methods which are completely analogous to those used in
\cite{OZ} for the investigation of a real differential calculus
on quantum Euclidean spaces.

In terms of the `angular momentum' generators
$M_j{^k}:=e_j{^k}-e_k{^j}$, the standard $r$-matrix on $\gotG =
so(n,\bR )$ is given by
\be
r =\frac{\er}{4} \sum_{j<k} (M_j{^k}+M_{j'}{^{k'}})\wedge
(M_j{^{k'}}+M_k{^{j'}})=
 \er \sum_{j<k<j'} M_j{^k}\wedge M_k{^{j'}},
\ee
where $j' := n+1 -j$ (the underlying Cartan subalgebra consists
of anti-diagonal matrices in this case).

 The above $r$-matrix is not quasi-triangular in the real sense
(this is a characteristic feature of compact simple groups).
Instead, one can find an invariant symmetric element $s$ of
$\gotG\otimes \gotG$ such that $w:= r \pm is$ satisfies the
classical Yang-Baxter equation.
We say that $(\gotG ,r)$ is {\em imaginary quasitriangular} in
this case.
In our case,
$$
s= \frac{\er}{2}\sum_{j,k} M_j{^k}\otimes M_k{^j}=\er \sum_{j,k}
(e_j{^k}\otimes e_k{^j} - e_j{^k}\otimes e_j{^k})
$$
(the simplest way to obtain $w$ is to extract the first order
term from the known $R$-matrix for the quantum $SO(n)$ group).
In order to satisfy (\ref{symsl}), we add $\lr I\otimes I$ to
$s$ with $\lr =\er$:
$$
s_\lr = \er \sum_{j,k}(e_j{^k}\otimes e_k{^j} -
e_j{^k}\otimes e_j{^k} + e_j{^j}\otimes e_k{^k}).
$$
It is clear that (\ref{general}) with $w$ replaced by $w_\lr := r
-is_\lr $ defines a complex-valued
Poisson structure on $T^*V\cong \bR ^{2n}$. Equally well we can treat
(\ref{general}) as defining a holomorphic Poisson structure on
the complexification $(T^*V)^{\mbC}\cong \bC ^{2n}$ of $T^*V$. In
the sequel we shall construct a real form of this holomorphic
Poisson manifold, playing the role of the phase space of
$(V,r_V)$.

We first derive Poisson brackets of coordinates with basic
$\gotG$-invariant functions:
$$
x^2 := g_{jk}x^jx^k,\qquad p^2 := p_jp_k g^{jk}, \qquad
E:=<p,x>=p_jx^j
$$
(summation convention assumed), where $g_{jk}$ is the metric
tensor (equal to the Kronecker delta in our orthonormal basis
$e_j$). Since $s$ is universal (the Killing form) ---
independent of $r$ (only the coefficient at $s$ depends on the
proportionality constant between $[r,r]$ and the canonical
element of $\Ltrzy \gotG$), the part
\be\label{univ}
\{ p_j,x^k \}_{\rm univ} = \dr _j {^k} +i\er (E\dr _j {^k} +
p_jx^k-x_jp^k )
\ee
of Poisson brackets (\ref{general}) corresponding to $\pi
_0-i(s_\lr )_{VV^*}$ is universal. It follows
that the brackets with $\gotG$-invariant functions are also
universal (these functions are Casimirs of $r_{T^*V}$ and it is
sufficient to use only $\pi _0 -i(s_\lr )_{VV^*}$). From
(\ref{univ}), it is now easy to obtain the following brackets:
\be\label{x2p2}
\{ x^2,x^j\} =0,\qquad \{ \frac12 p^2,x^j\} = p^j + i\er p^2
x^j,\qquad \{p^2 ,p^j\} =0 ,\qquad \{ p_j, \frac12 x^2\} = x_j
+i\er x^2 p_j,
\ee
$$
%\be\label{Exp}
\{ E,x^j\} = x^j +2i\er E x^j - i\er x^2 p^j,\qquad
\{ p_j,E\} = p_j +2i\er E p_j -i\er p^2 x_j.
%\ee
$$
Let us now introduce one more invariant function:
\be
\Lr := 1 +2i\er E - \er ^2 x^2p^2
\ee
(following the method of \cite{OZ}). From (\ref{x2p2})
 it follows that
\be
\{\Lr , x^j\} = 2i\er \Lr x^j,\qquad \{ p_j,\Lr \} = 2i\er
p_j\Lr .
\ee
Denote by $\Hol (Y)$ the algebra of holomorphic functions on a
complex manifold $Y$. Recall that any holomorphic map
$\phi\colon Y\to Z$ (of complex manifolds)
defines a linear multiplicative map $\Phi\colon \Hol (Z)\to \Hol
(Y)$ by the pullback: $\Phi (f) = f\circ \phi $. Similarly,
any anti-holomorphic map $\psi\colon Y\to Z$
defines an anti-linear multiplicative map
$\Psi\colon \Hol (Z)\to \Hol (Y)$ by  pullback followed by the
complex conjugation: $\Psi
(f) = \ov{f\circ \phi} $. Using $\Lr $, we define an
anti-holomorphic map $\psi $ from
$$
P^{\mbC}:= \{ (x,p)\in (T^*V)^{\mbC } : \Lr \neq 0\}\subset
(T^*V)^{\mbC }
$$
into $(T^*V)^{\mbC }$ by
\be
\Psi (x^j) = x^j,\qquad \Psi (p_j)=\frac{p_j +i\er p^2 x_j}{\Lr }.
\ee
Since
\be\label{Psi}
\Psi (p^2) = \frac{p^2}{\Lr },\qquad \Psi (E) = \frac{E+i\er
p^2x^2}{\Lr } ,\qquad \Psi (\Lr ) = \frac{1}{\Lr },
\ee
the underlying map $\psi $ maps $P^{\mbC}$ into $P^{\mbC}$.
Moreover, since $\Psi (\Psi (x^j))=x^j$ and
$$
\Psi (\Psi (p_j ))= \Psi \left(\frac{p_j +i\er p^2 x_j}{\Lr
}\right) = \Lr \left(\frac{p_j +i\er p^2 x_j}{\Lr } -i\er
\frac{p^2}{\Lr }x_j\right) = p_j,
$$
the anti-holomorphic map $\psi \colon P^{\mbC}\to P^{\mbC }$ is an
involution. Therefore we can define the corresponding {\em real
form} $P$ of $P^{\mbC}$ as the set of fixed points of $\psi$:
\be
P:= \{ z\in P^{\mbC} : \psi (z)=z\}.
\ee
The antilinear multiplicative involution $\Psi$ corresponding to
the map $\psi$ will be henceforth denoted by a star:
\be\label{star}
(x^j)^*=x^j, \qquad (p_j)^* = \frac{p_j +i\er p^2 x_j}{\Lr }.
\ee
Let us collect once more the basic formulas (\ref{Psi}):
$$
(p^2)^* = \frac{p^2}{\Lr },\qquad E^* = \frac{E+i\er
p^2x^2}{\Lr } ,\qquad  \Lr ^* = \frac{1}{\Lr },\qquad P:=\{
(x,p) : (x^j)^*=\ov{x^j},\;(p_j)^*=\ov{p_j} \}.
$$
The fundamental theorem of this section says that the star
operation is compatible with the Poisson brackets (\ref{general}).
\begin{Thm} The Poisson structure (\ref{general}) is real with
respect to the star operation (\ref{star}):
\be\label{reality}
\{ f,g\} ^* = \{ f^*,g^*\}.
\ee
\end{Thm}
\dowod
 We have to prove (\ref{reality}) in two cases: 1) $f=p_j$,
$g=x^k$ and 2) $f=p_j$, $g= p_k$. The case $f=x^j$, $g=x^k$ is
trivial.

1) We set $T_j := p_j +i\er p^2x_j$. Since
$$
\{p_j,x^k\} ^* = \dr _j{^k} - (p_l)^*(\ov{w_\lr})^{lk}_{jm}x^m
$$
and
$$
\{(p_j)^*,x^k\} =  \left\{ \frac{T_j}{\Lr },x^k\right\} =
\frac{1}{\Lr } \left( \{ T_j,x^k\} - 2i\er T_jx^k \right),
$$
we have to show that
\be\label{topr}
\{ T_j,x^k\} - 2i\er T_jx^k = \Lr \dr _j{^k} -
T_l(\ov{w_\lr})^{lk}_{jm}x^m .
\ee
Since
$$
\{ T_j,x^k\} = \dr _j{^k} - p_l (w_\lr )^{lk}_{jm}x^m + i\er
(2x_jp^k + 2i\er p^2 x_jx^k + p^2g_{jl}r^{lk}_{mn}x^mx^n),
$$
the left hand side of (\ref{topr}) equals
$$
 \dr _j{^k} - p_l (w_\lr )^{lk}_{jm}x^m + i\er
(2(x_jp^k -p_jx^k)  + p^2g_{jl}r^{lk}_{mn}x^mx^n),
$$
whereas the right hand side of (\ref{topr}) equals
$$
(1+2i\er E -\er ^2 x^2p^2)\dr _j{^k} - (p_l +i\er
p^2x_l)(\ov{w_\lr})^{lk}_{jm}x^m .
$$
Since
$$
p_l((w_\lr )^{lk}_{jm} - (\ov{w_\lr})^{lk}_{jm})x^m =
2p_l (-is_\lr )^{lk}_{jm}x^m = -2i\er (E\dr _j {^k} +
p_jx^k-x_jp^k )
$$
(cf. (\ref{univ})), it follows that (\ref{topr}) is equivalent to
$$
i\er p^2 g_{jl}r^{lk}_{mn}x^mx^n = -\er ^2 x^2 p^2\dr _j{^k}
-i\er p^2 x_l(\ov{w_\lr})^{lk}_{jm}x^m ,
$$
or,
$$
i g_{jl}r^{lk}_{mn}x^mx^n = -\er x^2\dr _j{^k}
-i x_l(\ov{w_\lr})^{lk}_{jm}x^m .
$$
Taking into account
\be\label{skew}
g_{jl}r^{lk}_{mn} = - g_{lm}r^{lk}_{jn}
\ee
($r^{lk}_{mn}$ belongs to $\gotG$ with respect to indices
$l,m$), it means that (\ref{topr}) is equivalent to
$$
ig_{lm}(is_\lr )^{lk}_{jn}x^mx^n = \er x^2 \dr _j{^k},
$$
which can be easily verified.

\vspace{1mm}

2) Since
$$
\{ T_j,\Lr \} = 2i\er T_j \Lr,
$$
we have
$$
\left\{ \frac{T_j}{\Lr} , \frac{T_k}{\Lr}\right\} =
\frac{\{T_j,T_k\}}{\Lr ^2}.
$$
Therefore it is sufficient to show that
$$
\{ T_j,T_k\} = T_lT_m r^{lm}_{jk}.
$$
We have
\begin{eqnarray*}
 \{ T_j,T_k\} & = & \{ p_j,p_k\} +i\er p^2 (\{p_j,x_k\} -\{
p_k,x_j\})+(i\er )^2\{ p^2x_j,p^2x_k\} \\
 & = & p_lp_m r^{lm}_{jk} + i\er p^2 [ -p_lr^{la}_{jm}x^mg_{ak} +
p_lr^{la}_{km}x^mg_{aj} +2i\er (p_jx_k-x_jp_k)] \\
 & & -\er ^2p^2
[p^2\{x_j,x_k\}+ 2( p_kx_j - x_kp_j ) ].
\end{eqnarray*}
Since
$$
-p_lr^{la}_{jm}x^mg_{ak} = p_lx_mr^{lm}_{jk}
$$
(by the argument similar to (\ref{skew})), we have finally
$$
 \{ T_j,T_k\} =p_lp_m r^{lm}_{jk} + i\er p^2 p_lx_m
(r^{lm}_{jk}-r^{lm}_{kj}) + (i\er p^2)^2 x_lx_mr^{lm}_{jk}
= (p_l +i\er p^2x_l)(p_m+i\er p^2 x_m)r^{lm}_{jk}
{}.
$$

\qed

{\em Corollary:} \ $P$ is endowed with a structure of a real
analytic Poisson manifold. If $f_0,g_0$ are two real analytic
functions on $P$, then their Poisson bracket is defined by
$$
\{ f_0 , g_0 \} := \{ f,g\}| _P ,
$$
where $f,g$ are the (local) holomorphic extensions of $f_0,g_0$
to $P^{\mbC }$ (by (\ref{reality}), the restriction of $\{ f,g\}$
to $P$ is real).

$P$ is the required phase space of $(V,r_V )$.

\section{Poisson action of $SU(n)$ on $\dbC ^n$}

We treat here $V=\bC ^n$ as a real manifold ($V\cong
\bR^{2n}$). Specifying (\ref{stand}) to the case of $SU(n)$ we
get the following standard $r$-matrix (see (\ref{sunbas}) for
the basis of $su (n)$)
\be\label{sunstan}
r =\er \frac12 \sum_{j<k} (e_j{^k}-e_k{^j})\wedge
J(e_j{^k}+e_k{^j}) ,
\ee
where $J\colon V\to V$ is the complex structure of $V$
(multiplication by the imaginary unit).
{}From now on we set $\er =1$ (arbitrary
$\er$ will be restored in the final formulas). It is convenient
to work with the complexification $V^{\mbC}\cong V\oplus iV$ and
the complex-linear embedding
$$
     V\ni z\mapsto z^{\mbC} := \frac12 (z -iJz)\in V^{\mbC} .
$$
We have
$$
 z = z^{\mbC} + \ov{z^{\mbC}}, \qquad Jz = i(z^{\mbC} -\ov{z^{\mbC}})
$$
and the typical notation:
$$
(e_k)^{\mbC} = (\frac{\partial}{\partial x_k})^{\mbC} =
\frac{\partial}{\partial z_k} = \partial _k.
$$
Note that
$$
e_j{^k}z = (e_j{^k}z)^{\mbC} + \ov{(e_j{^k}z)^{\mbC}}=
           (e_j z^k)^{\mbC} + \ov{(e_jz^k)^{\mbC}}
	   = z^k\partial _j + \cz ^k\cpartial _j.
$$
In this notation, the fundamental bi-vector field $r_V$ looks as
follows
\be
r_V (z) = i\sum_{jk} \sgn (k-j) \left( \frac12 \nabla
_j\wedge\nabla _k -  \frac12 \cnabla_j\wedge\cnabla _k	 +
|z^j|^2 \partial _k\wedge \cpartial _k\right) ,
\ee
where $\nabla _k := z^k\partial _k$, $\cnabla _k := z^k\cpartial
_k$.
\begin{Lem}
\ $[r_V,r_V] (z) = -||z||^2 Jz\wedge \pi _0$,\ \ where
$$
 \pi _0 = 2i\sum_k \partial _k\wedge \cpartial _k
$$
 is the canonical constant bi-vector on $V =\bC ^n =\bR ^{2n}$.
 \end{Lem}
\dowod
Taking into account
$$ [|z^j|^2\partial _k\wedge \cpartial _k, \nabla _a\wedge
\nabla _b ] =  \cz ^j (-\cpartial _k)\wedge [z^j\nabla _k,\nabla
_a\wedge \nabla _b]
$$
$$
=-\cz ^j \cpartial _k\wedge [(z^j\dr _k^a\partial _a - z^a\dr
_a^j\partial _k)\wedge\nabla_b + \nabla _a\wedge (z^j\dr
_k^b\partial _b - z^b\dr _b^j\partial _k)]
$$
$$
=-|z ^j|^2 \cpartial _k\wedge\partial _k\wedge
[(\dr _k^a -\dr _a^j)\nabla _b - (\dr _k^b -\dr _b^j)\nabla _a ]
$$
and
$$
 [|z^j|^2\partial _k\wedge \cpartial _k, |z^a|^2\partial
_b\wedge \cpartial _b] = [z^j\partial _k\wedge \cz ^j\cpartial
_k , z^a\partial _b\wedge \cz ^a\cpartial _b]=
$$
$$
= [z^j\partial _k, z^a\partial _b]\wedge
 \cz ^j\cpartial _k\wedge \cz ^a\cpartial _b
+ z^j\partial _k\wedge z^a\partial _b\wedge
[ \cz ^j\cpartial _k, \cz ^a\cpartial _b]=
z^j\partial _k\wedge z^a\partial _b\wedge
( \cz ^j\dr _k^a\cpartial _b- \cz ^a\dr _b^j\cpartial _k)+ c.c.
$$
($+ c.c.$ means `plus complex conjugated terms'), we see that
$[r_V,r_V](z)$ equals
$$
-\sum_{jkab}\sgn (k-j)\sgn (b-a)
\{ |z^j|^2[(\dr _k^a -\dr _a^j)\nabla _b - (\dr _k^b -\dr
_b^j)\nabla _a ] + 2 |z^a|^2 \dr _b^j \nabla _b\} \wedge
\partial _k\wedge \cpartial _k \, + \, c.c. =
$$
$$
= -2 \sum_{jkab}\sgn (k-j)\sgn (b-a)
[ |z^j|^2(\dr _k^a -\dr _a^j)\nabla _b
+  |z^a|^2 \dr _b^j \nabla _b] \wedge
\partial _k\wedge \cpartial _k \, + \, c.c.
$$
$$
= -2\sum_{kb} \sum_{ja}\sgn (k-j)\sgn (b-a)
[ |z^j|^2(\dr _k^a -\dr _a^j)
+  |z^a|^2 \dr _b^j] \nabla _b \wedge
\partial _k\wedge \cpartial _k \, + \, c.c.
$$
Note that
$$
\sum_{ja}\sgn (k-j)\sgn (b-a)
( |z^j|^2\dr _k^a -|z^j|^2\dr _a^j
+  |z^a|^2 \dr _b^j) \nabla _b \wedge
\partial _k\wedge \cpartial _k
$$
$$
= \sum_j |z^j|^2 (\sgn (k-j)\sgn (b-k) + \sgn (b-k)\sgn (j-b) +
\sgn (j-b)\sgn (k-j))
$$
and
$$
\sgn (k-j)\sgn (b-k) + \sgn (b-k)\sgn (j-b) +
\sgn (j-b)\sgn (k-j) = -1
$$
for $b\neq k$. It follows that
$$
[r_V, r_V] (z) = 2 ||z||^2 \sum_b (\nabla _b - \cnabla _b)\wedge
\sum_k \partial _k\wedge \cpartial _k = -||z||^2
\sum_b i(\nabla _b - \cnabla _b)\wedge
\sum_k 2i\partial _k\wedge \cpartial _k.
$$

\qed

{\em Corollary:} \ For any classical $r$-matrix $\tilde{r}$ on
$\gotG =su(n)$ there is a constant $c$ such that
$[\tilde{r}_V,\tilde{r}_V]= -c||z||^2 Jz\wedge \pi_0$. Poisson
structures $\pi _V$ on $V$ for which the action of $SU(n)$ on $V$
is Poisson are exactly bi-vector fields
\be\label{Poi}
\pi _V = \tilde{r}_V + \Dr,
\ee
such that the bi-vector field $\Dr$ on V is $SU(n)$-invariant
and satisfies
\be\label{braDr}
[\Dr ,\Dr ](z) = c||z||^2 Jz\wedge \pi_0 .
\ee

It is easy to show that all $SU(n)$-invariant bi-vector fields
$\Dr $ on $V$ are of the following form:
\be\label{form}
\Dr  = \frac12 a\pi _0 + \frac12 bz\wedge Jz,
\ee
where $a=a(||z||^2)$, $b=b(||z||^2)$ are arbitrary functions of
$||z||^2$. We shall write condition (\ref{braDr}) in terms of
these functions.
\begin{Lem}
\ $[\Dr ,\Dr ](z) = ||z||^2 Jz\wedge \pi_0$ \ if and only if
\be\label{warunek}
 aa' + b(a - a't) =t.
 \ee
 Here $t\equiv ||z||^2$ and prim means the differentiating with
respect to the variable $t$.
\end{Lem}
\dowod
 If $K$, $L$ are bi-vector fields and $f$, $g$ are functions, then
\be
[fK,gL] = fg[K,L] -fK(dg)\wedge - gK\wedge L(dg),
\ee
where by $K(dg)$ we denote the contraction of $K$ with $dg$ on
the first place. In particular,
$$
[fK,fK] = f^2 - 2fK\wedge K(df).
$$
Using
$$
 \pi_0 \left(\frac12 d ||z|| ^2\right) = -Jz,\qquad (z\wedge Jz)
\left(\frac12 d ||z|| ^2\right)= ||z||^2Jz,
$$
and
$$
[\pi _0 , z\wedge Jz ] = 2Jz\wedge \pi _0,
$$
we obtain
\begin{eqnarray*}
 {}[ \frac12 a \pi _0 , \frac12 a \pi _0 ] & = & aa' Jz\wedge \pi _0\\
 {}[bz\wedge Jz,bz\wedge Jz ] & = & 0 \\
 {}[\frac12 a \pi _0 ,\frac12 b z\wedge Jz ] & = & \frac12 b(a
-a'||z||^2)Jz\wedge \pi _0 .
\end{eqnarray*}
{}From this, (\ref{warunek})  follows immediately.

\qed

Of course, the easy way to solve (\ref{warunek}) is to write
\be\label{solution}
b=\frac{t-aa'}{a-a't},
\ee
but in such a way we have no control of regularity of those
functions and we do not see the simplest cases.
To pick up the simplest cases, let us consider $\Dr $ at most
quadratic, i.e.
$a=a_0+a_1t$, $b=b_0$, where $a_0,a_1,b_0$ are some constants.
Inserting this form of $a$ and $b$ in (\ref{warunek}) gives the
following two cases:

1. \ $a_0 = 0$, $a_1 =\pm 1$, $b_0$ arbitrary,

2. \ $a_0$ arbitrary, $a_1=\pm 1$, $b_0=\mp 1$.

One of the simplest non-quadratic solution for $\Dr $ is the
following solution of degree 4:

3. \ $a = h=const \neq 0$, $b = \frac{1}{h} t$

An example of a non-singular rational solution is given by
$a= 1-t^2 $ (in this case the denominator of (\ref{solution}) is
positive: $a-a't=1+t^2$).

Another way to pick up a simple case is to assume that $\Dr $ is
(as $r_V$) tangent to spheres $||z||=const$. It is easy to show
that this conditions holds if and only if $a=bt$. In this case
(\ref{warunek}) reduces to
$$ab=t.$$
It means that $a=\pm t$, $b=\pm 1$. This is a special case of type
1 above.

We end by listing the explicit form of the Poisson brackets
corresponding to the mentioned cases. From the general form
(\ref{Poi}), (\ref{form}) with the standard $r$-matrix (\ref{sunstan}),
we obtain
\begin{eqnarray}
\{ z^j,z^k\} & = & i\er z^jz^k \qquad \mbox{for}\;\; j<k \\
\{ z^j,\cz ^k\} & = & -i\er b z^j\cz ^k \qquad \mbox{for}\;\;
j\neq k \\
\{ z^j, \cz ^j\} & = & i\er \sum_k \sgn (j-k) \cdot |z^k| \, +
\, i\er a - i\er b |z^j|
\end{eqnarray}
(we have restored the parameter $\er$). Now we list the cases
which seem to be most intersting ones.
\begin{enumerate}
\item {\em Poisson $SU(n)$-spheres}. \ According to the
discussion above, there are only two Poisson structures on $\bC
^n$ solving our problem and tangent to spheres, namely
\begin{eqnarray*}
\{ z^j,z^k\} & = & i\er z^jz^k \qquad \mbox{for}\;\; j<k \\
\{ z^j,\cz ^k\} & = & -i\sr\er b z^j\cz ^k \qquad \mbox{for}\;\;
j\neq k \\
\{ z^j, \cz ^j\} & = & i\er (\sr ||z||^2 -\sr |z|^j + \sum_{k}
\sgn (j-k) \cdot |z^k|) = 2\sr i\er\sum_{\sr k <\sr j} |z^k|^2 ,
\end{eqnarray*}
where $\sr = \pm 1$. The function $z\mapsto ||z||^2$ is a
Casimir function of this Poisson structure (and can be fixed,
which leads to a sphere $S^{2n-1}$).
\item
{\em Twisted annihilation and creation `operators'}. \
Setting $h=\er a_0$ in the case 2 above, we obtain
\begin{eqnarray*}
\{ z^j,z^k\} & = & i\er z^jz^k \qquad \mbox{for}\;\; j<k \\
\{ z^j,\cz ^k\} & = & i\sr\er b z^j\cz ^k \qquad \mbox{for}\;\;
j\neq k \\
\{ z^j, \cz ^j\} & = & ih +i\er (\sr ||z||^2 +\sr |z|^j + \sum_{k}
\sgn (j-k) \cdot |z^k|) = ih + 2\sr i\er\sum_{\sr k \leq \sr j}
|z^k|^2 ,
\end{eqnarray*}
where $\sr = \pm 1$. This is the Poisson version of the
`twisted canonical commutation relations' of \cite{PuW} (see
also (\cite{Sh})). It may describe the phase space of a Poisson deformed
harmonic oscillator.
\item The non-quadratic brackets corresponding to the case 3
above are given by
\begin{eqnarray}
\{ z^j,z^k\} & = & i\er z^jz^k \qquad \mbox{for}\;\; j<k \\
\{ z^j,\cz ^k\} & = & -i\frac{\er}{h} ||z||^2 z^j\cz ^k \qquad
\mbox{for}\;\; j\neq k \label{deg4}\\
\{ z^j, \cz ^j\} & = & i\er h +i\er (-\frac{1}{h} ||z||^2 |z^j|^2
+\sum_{k}\sgn (j-k) \cdot |z^k|) .\label{deg41}
\end{eqnarray}
\end{enumerate}

{\bf Problem:} \ What is the quantum counterpart of condition
(\ref{warunek})? The quantum counterpart of relations
(\ref{deg4}), (\ref{deg41})?

\end{document}